\documentclass[3p,times,reqno]{elsarticle}

\usepackage{natbib,url,hyperref,doi}
\usepackage{ecrc}

\usepackage{enumitem}

\volume{00}

\firstpage{1}

\journalname{Computers \& Mathematics with Applications}

\runauth{R. M. Evans et al.}

\jid{camwa}

\jnltitlelogo{\small{Computer \\ \& \\ Mathematics \\with Applications}}

\usepackage{amssymb}
 \usepackage{amsthm}

\usepackage{xcolor,soul} 
\definecolor{blue1}{RGB}{66, 226, 244}

\usepackage{graphicx,subfig,epstopdf}  
\usepackage{epsfig,epstopdf,verbatim,subfig,caption}
\usepackage[english]{babel}
\usepackage[export]{adjustbox}%
\usepackage[utf8]{inputenc}
\usepackage{amsmath}
\usepackage{graphicx}
\usepackage[colorinlistoftodos]{todonotes}
\usepackage{lmodern}
\usepackage{amssymb,amsthm}
\usepackage{thmtools}
\usepackage{verbatim}
\usepackage{filecontents}

\usepackage{chemarr}
\usepackage{textcomp}  

\newcommand{\Da}{\text{Da}}

\newcommand{\p}{\partial}
\newcommand{\wt}{\widetilde}

\usepackage{epsfig}
\usepackage{psfrag}
\usepackage{float}
\usepackage{caption}

 \biboptions{sort&compress}
\newtheorem{theorem}{Theorem}[section]
\newtheorem{lemma}[theorem]{Lemma}

\theoremstyle{definition}
\newtheorem{definition}[theorem]{Definition}

\newtheorem{example}[theorem]{Example}

\newtheorem*{acknowledgement}{Acknowledgement}

\theoremstyle{remark}
\newtheorem{remark}[theorem]{Remark}

\begin{document}
\begin{frontmatter}
 \title{Applications of fractional calculus in solving Abel-type integral equations: Surface-volume reaction problem}  
\author{Ryan M. Evans\tnoteref{fn1}}
\tnotetext[fn1]{Current address: Applied and Computational Mathematics Division, Information and Technology Lab, National Institute for Standards and Technology, Gaithersburg, MD 20899, USA}
%
\author{Udita N. Katugampola\corref{corresp}}
%
\author{David A. Edwards}
%

\cortext[corresp]{Corresponding author.  Tel.: +13028312694.\\
\indent {\it E-mail address:} uditanalin@yahoo.com (U. N. Katugampola)}

\address{Department of Mathematical Sciences, University of Delaware, Newark, DE 19716, U.S.A.}

\begin{abstract}
In this paper we consider a class of partial integro-differential equations of fractional order, motivated by an equation which arises
as a result of modeling surface-volume reactions in optical biosensors.  We solve these equations by employing techniques from fractional calculus;  several examples are discussed.  
Furthermore, for the first time, we encounter an order of the fractional derivative other than $\frac{1}{2}$ in an applied problem. 
Hence, in this paper we explore the applicability of fractional calculus in real-world applications, further strengthening the true nature of fractional calculus.
\end{abstract}
\begin{keyword}
Abel integral equation \sep Riemann--Liouville integral \sep Caputo fractional derivative \sep Fractional calculus \sep Surface-volume reactions 
%
%
\MSC[2010] 26A33 \sep 92C45
\end{keyword}
\end{frontmatter}
%

%
\section{Introduction}
%
%

Until very recently, the fractional calculus had been a purely mathematical tool without apparent applications. {Currently, fractional 
dynamical equations play a major role in modeling of anomalous behavior and memory effects, which are common characteristics of 
natural phenomena \cite{Herr,what,quo}. The fact that fractional derivatives introduce a convolution integral with a power-law memory kernel makes the 
fractional differential equations an important model to describe memory effects in complex systems. Thus, it is seen that fractional derivatives or integrals appear
 naturally when modeling long-term behaviors, especially in the areas of viscoelastic materials and viscous fluid dynamics \cite{kumar,yang}.

Abel's study of the tautochrone problem  \cite{Abel} is considered to be the first application of fractional calculus to an engineering problem.  In it one finds the path where the time it takes for an object to fall under the influence of gravity is independent of the initial position. The solution, which was solved using a fractional calculus approach, is now known to be a part of the inverted cycloid \cite{Abel,Ross}.

Now it is not hard to find very interesting and novel applications of fractional differential equations in physics, chemistry, biology, engineering, finance and other areas of sciences that have been developed in the last few decades. Some of the applications include: diffusion processes \cite{mainardi1,sie}, mechanics of materials \cite{tor,Caputo}, combinatorics \cite{udita3,oeis}, inequalities \cite{Chen1}, analysis \cite{u-1}, calculus of variations \cite{mal,u-3,u-4,u-5,u-5-2,Almeida1-doi}, signal processing \cite{mark},
image processing  \cite{bai}, advection and dispersion of solutes in porous or fractured media \cite{ben}, modeling of viscoelastic materials under external forces \cite{fre}, bioengineering \cite{mag}, relaxation and reaction kinetics of polymers \cite{gloc}, random walks \cite{gro1}, mathematical finance \cite{gro2}, modeling of combustion \cite{led}, control theory \cite{podl}, heat propagation \cite{podl2}, modeling of viscoelastic materials \cite{jar} and even in areas such as psychology \cite{ahmad,song}. The list is by no means complete. It is easy to find hundreds, if not thousands, of new applications in which the fractional calculus approach is more than welcome. 

This is the first in a series of papers that seeks to find further potential applications of fractional calculus in solving real-world problems, a journey that can benefit both the understanding of profound complexities in the application, and the field of fractional calculus itself. As an application of the theory developed in this paper, we consider the surface-volume reaction problem. The governing equations of the mathematical formulation of such models naturally give rise to a nonlinear equation that contains a fractional integral embedded in it, and which has no solutions to date. Thus, in this paper we both extend the theory of fractional calculus methods by considering equations motivated by modeling the surface-volume reactions, and explore another interpretation of the fractional integral.

The remainder of this paper is organized as follows. The definitions and basic results are given in Section~\ref{Section: Basic Definitions and Preliminary Results}. 
In Section~\ref{Section: Main Results}, we give the main results, which are generalizations of Abel's integral approach to the tautochrone problem. In 
Section~\ref{Section: Illustrative Examples}, we give illustrative examples to motivate our approaches. 
One of the main examples is the surface-volume reaction problem that has several very interesting applications in mathematical biology and engineering \cite{Edwards2}.
 \section{Basic Definitions and Preliminary Results}\label{Section: Basic Definitions and Preliminary Results}

We adopt definitions given in \cite{pod} or in the encyclopaedic book by Samko {\it et al}. \cite{samko} here.  We begin by introducing the concept of a {\it Riemann--Liouville fractional integral}:

\begin{definition}[\cite{pod}]Let $\alpha >0$ with $n-1 < \alpha \leq n, \; n \in \mathbb{N}$, and $a < x<b$. The left- and right- Riemann--Liouville fractional integrals of order $\alpha$ of a function $f$ are given by 
\[
    J^\alpha_{a+}f(x) = \frac{1}{\Gamma(\alpha)}\int_{a}^x (x-t)^{\alpha -1} f(t) \, dt \quad \mbox{and} \quad J^\alpha_{b-}f(x) = \frac{1}{\Gamma(\alpha)}\int_x^b (t-x)^{\alpha -1} f(t) \, dt 
\]
\noindent respectively, where $\Gamma(\cdot)$ is Euler's gamma function defined by
$$
\Gamma(x) = \int_0^\infty t^{x-1}e^{-t} \, dt.
$$
\end{definition}

Non-local fractional derivatives are defined via fractional integrals \cite{samko,udita1,udita2}, while the local fractional derivatives are defined via a limit-based approach \cite{kol1,kol2}. A new class of controlled-derivative approach appeared in \cite{katu}. A criteria to test whether a given derivative is a fractional derivative appeared in \cite{wfra,udita4}. Among other approaches, in this work, we utilize only the non-local Volterra-type definitions for the fractional derivative given below.

\begin{definition}[\cite{pod}] The left- and right- {\it Riemann--Liouville fractional derivatives} of order $\alpha > 0$, $n-1 < \alpha < n$, $n \in \mathbb{N}$, are defined by
\begin{eqnarray*}
    D^\alpha_{a+}f(x) &=& \frac{1}{\Gamma(n-\alpha)} \bigg(\frac{d}{dx}\bigg)^n \int_{\rm a}^x (x-t)^{n-\alpha -1} f(t) \, dt,\\
    D^\alpha_{b-}f(x) &=& \frac{(-1)^n}{\Gamma(n-\alpha)}\bigg(\frac{d}{dx}\bigg)^n\int_x^b (t-x)^{n-\alpha -1} f(t) \, dt,
\end{eqnarray*}
\end{definition}
\noindent
respectively. It can be shown that in the case of $\alpha \in \mathbb{N}$ the above definitions coincide with the standard definition of the $n^{th}$-derivative of $f(x)$. 

\begin{definition}[\cite{pod}] The left- and right- {\it Caputo fractional derivatives} of order $\alpha > 0, n-1 < \alpha < n, \, n \in \mathbb{N}$, are defined by
\[
    {}^{\rm C} D^\alpha_{a+}f(x) = \frac{1}{\Gamma(n-\alpha)} \int_{\rm a}^x (x-t)^{n-\alpha -1} f^{(n)}(t) \, dt \quad \mbox{and} \quad {}^{\rm C} D^\alpha_{b-}f(x) = \frac{(-1)^n}{\Gamma(n-\alpha)}\int_x^b (t-x)^{n-\alpha -1} f^{(n)}(t) \, dt,
\]
\end{definition}
\noindent
respectively. It can be shown that in the case of $\alpha \in \mathbb{N}$ the above definitions reduce to the standard definition of the $n^{th}$-derivative of $f(x)$. To see this, let us assume that $0 \leq n-1 < \alpha < n$, and $f(x)\in C^{n+1}[a, T]$. Then in the case of Caputo's derivative, we have, by integration by parts \cite[p.~79]{pod}, 
\begin{align}
\lim_{\alpha \rightarrow n}  {}^{\rm C} D^\alpha_{a+}f(x) &= \lim_{\alpha \rightarrow n} \left[\frac{f^{(n)}(a)(x-a)^{n-a}}{\Gamma(n-\alpha+1)} + \frac{1}{\Gamma(n-\alpha+1)}\int_{\rm a}^x(x-\tau)^{n-\alpha}f^{(n+1)}(\tau)\,d\tau\right] \\
             &= f^{(n)}(a) + \int_{\rm a}^x f^{(n+1)}(\tau)\,d\tau = f^{(n)}(x), \quad n=1, 2, \ldots.
\end{align}
This shows that the Caputo derivative is a generalization of the integer-order derivative. A different proof of the fact in question, which does not use integration by parts, can be found in \cite[pp.~49, 51]{kai} using an equivalent definition of the Caputo derivative. On the other hand, for the Riemann--Liouville derivative, the result follows easily since $J^0_{a+}$ is defined to be the identity operator and one can use the fact that $D^\alpha_{a+} = D^{\lceil\alpha\rceil}_{a+}J^{\lceil\alpha\rceil -\alpha}_{a+}$ \cite[p.~27]{kai}.

The following relations between the Caputo and Riemann--Liouville fractional derivatives are noteworthy \cite{pod,kai}:

\begin{equation}
{}^{\rm C} D^\alpha_{a+} f(x)= D^\alpha_{a+} f(x) - \sum^{n-1}_{k=0} \frac{f^{(k)}(a)}{\Gamma(k-\alpha+1)} (x-a)^{k-\alpha}
\end{equation}
and
\begin{equation}
{}^{\rm C} D^\alpha_{b-} f(x)= D^\alpha_{b-} f(x) - \sum^{n-1}_{k=0} \frac{f^{(k)}(b)}{\Gamma(k-\alpha+1)} (b-x)^{k-\alpha}.
\end{equation}
The two relations in question have interesting properties. That is, if $f \in C^n[a, b]$ and $f^{(k)}(a) = 0$, $ k = 0, 1, \ldots, n-1,$ then
\[
   {}^{\rm C} D^\alpha_{a+} f = D^\alpha_{a+} f,
\]
and if $f^{(k)}(b) = 0$, $ k = 0, 1, \ldots, n-1,$ then
\[
   {}^{\rm C} D^\alpha_{b-} f = D^\alpha_{b-} f.
\] 
The fractional integrals and derivatives also satisfy the following important properties: fractional operators are linear, that is, if $L$ is a fractional integral or derivative, then
\[
    L(f + kg) = L(f) + kL(g)
\]
for any functions $f, g \in C^n[a, b]$ or $f, g \in L^p(a, b)$ (as the case may be) and $k \in \mathbb{R}$. For any $\alpha, \beta >0$, they also satisfy the following semigroup properties:
\[
     J^\alpha J^\beta = J^{\alpha +\beta} \qquad \mbox{and} \qquad D^\alpha D^\beta = D^{\alpha +\beta}.
\]
Further, if $f \in L^\infty (a, b)$ or $f \in C^n[a, b]$ and $\alpha > 0$, then \cite[p. 95]{samko}
\begin{equation}
   {}^{\rm C} D_{a+}^\alpha J^\alpha_{a+} f = f \qquad \mbox{and} \qquad {}^{\rm C} D_{b-}^\alpha J^\alpha_{b-} f = f.
\label{prop1}
\end{equation}
Equation \eqref{prop1} will be the key property used to prove the main results of this paper. It also says that the Caputo derivative is the left-inverse of the Riemann--Liouville fractional integral. Unfortunately, the Caputo derivative is not the right-inverse of the Riemann--Liouville integral. But, we have the following result concerning those two operators \cite{kai,pod}: 
\[
   J^\alpha_{a+}{}^{\rm C} D_{a+}^\alpha f(x) = f(x) - \sum^{n-1}_{k =0} \frac{f^{(k)}(a)}{k!}(x-a)^k,
\]
and
\[
 J^\alpha_{b-}{}^{\rm C} D_{b-}^\alpha f(x) = f(x) - \sum^{n-1}_{k =0} \frac{f^{(k)}(b)}{k!}(b-x)^k,
\]
if $f \in C^n[a, b]$ and $\alpha > 0$. 



The Caputo derivatives have some advantage over the Riemann--Liouville derivatives. One of the most important properties is that the Caputo derivative of a constant function is zero whereas the Riemann--Liouville derivative is not. From an analytical point of view, this result does not create difficulties though it is commonly not acceptable in the physical sense. 

Using a convolution integral approach, Podlubny \cite{pod1,pod2} provided geometric and physical interpretations of fractional integration and fractional differentiation for the Riemann--Liouville fractional operators, Caputo derivatives, Riesz potentials and Feller potential. 
Several other attempts to provide geometric meanings to fractional operators can be found, for example, in the works by Ben Adda \cite{adda1,adda2}. Further approaches can be found in,
for example, \cite{Ren,Gor1,Hey,Kirk1,Mon,Nig,Rut,Rut2,Yu1,Mac,recent1,recent2,herr2,Main2}. Among such approaches, as pointed out by Rutman \cite{Rut}, there were misunderstandings in some
 approaches which try to find interpretations of fractional integrals and derivatives \cite{Nig}. Thus, one of the goals of this paper is to find concrete examples that can be directly 
related to the fractional operators so that we may find a better physical interpretation rather than building our own models. Further details of the approach can be found in
 Section \ref{Section: Illustrative Examples} of this work. 

The rest of the paper will adopt the following definitions of Volterra-type integrals.
\begin{definition}[\cite{jake}]
An integral equation of the form
\begin{equation}\label{eq1}
w(t) = \int_{\rm a}^t K(t, s)v(s)\,ds, \quad a \leq s \leq t \leq T,
\end{equation}  
is known as a \textit{Volterra integral equation of the first kind} with kernel $K(t, s)$. Its kernel is said to be \textit{non-singular} or \textit{singular} depending on whether $K(t, s)$ is continuous or discontinuous on the triangular region $a \leq s \leq t \leq T$. In the special case where $K(t, s)$ has the form 
\[
   K(t, s) = \frac{k(t, s)}{(t - s)^\alpha}, \quad 0 <\alpha <1,
\]
with $k(t, s)$ continuous on $a \leq s \leq t \leq T$, (\ref{eq1}) is called a \textit{generalized Abel integral equation} or an \textit{integral equation of Abel type}. 
\end{definition}

For Abel type integral equations of the form
\begin{equation}\label{eq2abel}
  w(t) = \int_0^t \frac{k(t, s)}{(t - s)^\alpha}v(s) \, ds, \quad 0 \leq s \leq t \leq T, \; 0 < \alpha < 1,
\end{equation}
we have the following result \cite[pp.~80--82]{kow1}:
\begin{theorem}\label{thm1} Assume $w(t)$ is defined as in Eq.\ \eqref{eq2abel}.  If
	\begin{enumerate}
		\item [$(i)$] $k(t, t) \ne 0, \quad 0 \leq t \leq T$,
		
		\item [$(ii)$] $k(t, s)$ has continuous partial derivatives up to order $m \in \{0, 1, 2, \ldots \}$, and $\partial^{m+1}k(t, s)/\partial t^{m+1}$ is continuous on $0 \leq s \leq t \leq T$,
		\item [$(iii)$] $w(t) \in C^{m+1}[0, T]$ and $w^{(l)}(0) = 0, \, (l = 0, 1, 2, \ldots m)$, and
		
		\item [$(iv)$] $W(t) = \frac{d}{dt}\int_0^t \frac{w(s)}{(t - s)^{1-\alpha}}\, ds \in C[0, T]$,
	\end{enumerate}
then $v\in C^m[0, T]$ is unique.  
\end{theorem}

It is interesting to see that condition (iv) of Theorem \ref{thm1} is the Riemann--Liouville fractional derivative of order $\alpha$ for $0 < \alpha < 1$ though there is no explicit mention about such a derivative in the literature.  

Now consider the following Abel integral equation of the first kind:
\begin{equation}
     f(x) = \int_0^x \frac{g(t)}{(x-t)^\alpha}\,dt, \quad 0<\alpha<1, \quad 0 \leq x \leq b, \label{prop2}
\end{equation}
where $f \in C^1[a, b]$ is a given function satisfying $f(0) = 0$ and $g$ is the unknown function to be determined.
\noindent Several methods for solving such Volterra-type integral equations include: the Chebyshev polynomial approach \cite{avaz}, orthogonal polynomial approach \cite{Mine}, Galerkin methods \cite{paul}, collocation methods \cite{brun}, Haar/CAS wavelet method \cite{lep,sae,sae2} and Mikusinski's 
operator approach \cite{Li}. Several classical approaches in solving integral equations can be found in the book by Kanwal \cite{Kanwal}. 

Jahanshahi {\it et al}. \cite{Salman} solved the problem given in \eqref{prop2} using a fractional calculus approach and arrived at the solution given below. 
\begin{theorem}[\cite{Salman}] The solution to \eqref{prop2} is 
\[
    g(x) = \frac{\sin (\pi\alpha)}{\pi} \int_0^x \frac{f^\prime (t)}{(x-t)^{1-\alpha}}\,dt.
\]
\end{theorem}
The interested reader may find a similar result in \cite{Li} and also in the work by Gelfand and Vilenkin \cite[Section 5.5]{Gal}. Those approaches are very similar to the method introduced by Abel in \cite{Abel}. 
Continuing in the same direction, we shall solve a Volterra-type integral equation that involves two variable functions in the integrand, which contains partial derivatives. The method we describe here is new, to the authors' knowledge. 

Before elaborating on the new results, let us make some remarks about the notations that will be used in the rest of the paper. 

\subsection{Notations}
Here we will go through some notation that is common to the current literature of the subject.  It will be convenient to use operator notation. By $J^\alpha$ we will mean the fractional integral operator
\[
J^\alpha f(x,t)=\frac{1}{\Gamma(\alpha)}\int_0^x\frac{f(\nu,t)}{(x-\nu)^{1-\alpha}}\, d \nu,\quad \alpha\in [0,1].
\]
Note this fractional integral is with respect to space.  If we take the fractional integral with respect to time, we will denote this by
\[
J_t^\alpha f(x,t)=\frac{1}{\Gamma(\alpha)}\int_0^t\frac{f(x,\tau)}{(t-\tau)^{1-\alpha}}\, d \tau.
\]
Also by $If(x,t)$, or simply by $If$, we will always mean an integral with respect to time:
\[
If=\int_0^tf(x,\tau)\, d\tau.
\]  
Additionally we will use the notation
\[
Igf=\int_0^tg(\tau)f(x,\tau)\, d\tau.
\]
When dealing with fractional derivatives, we will always understand the fractional derivative in the Caputo sense
\[
D^\alpha f=\frac{1}{\Gamma(1-\alpha)}\int_0^x\frac{1}{(x-\nu)^{\alpha}}\frac{\p f(\nu,t)}{\p x}\, d \nu.
\]
We use a similar definition for the temporal fractional Caputo derivative
\[
D^\alpha_t f=\frac{1}{\Gamma(1-\alpha)}\int_0^t\frac{1}{(t-\tau)^{\alpha}}\frac{\p f(x,\tau)}{\p t}\, d \tau.
\]

We are now ready to state our results.
\section{Main Results}\label{Section: Main Results}
The following is the main result of this paper, which is motivated by a surface-volume reaction problem after taking into account transport effects. 
There equations are nonlinear in nature and have no exact solutions to date. Equations of this form can be found in \cite{Edwards1,Edwards2,Edwards3}. 
For example, the governing equations of this form can be found in Equations (2a) and (2b) of \cite{Edwards3}. 
Further treatment of this class of equations can be found in Section \ref{Setting} of this work. 

\begin{lemma}\label{lemma1}
Let $C$ be a function, which is a linear combination of separable functions of $x$ and $t$ and functions of only $x$ (but, not constant functions) such that its first partial derivative with respect to $x$ exists. Let $0 <\alpha <1$. Then the solution of the integral equation
\begin{equation}\label{eq2}
   C(x, t) = \int_0^x \frac{\partial B}{\partial t}(x-\xi, t)\, \frac{d\xi}{\xi^\alpha}, \quad \mbox{with} \quad B(x, 0) = 0,	
\end{equation}
is given by
\begin{equation}\label{eq3}
   B(x, t) = \frac{\sin (\pi \alpha)}{\pi}\int_0^t \int_0^x \frac{\partial C}{\partial x}(\xi, \tau) \frac{d\xi}{(x-\xi)^{1-\alpha}}\,d\tau.
\end{equation}
\end{lemma}
\begin{proof}
We prove this first using a change of variable and then identifying (\ref{eq2}) as a fractional integral of order $\alpha$.
 To that end, let $x-\xi= \eta$ and $F(x;t)=\frac{\partial B}{\partial t}(x, t)$. Making this substitution, (\ref{eq2}) becomes
\begin{equation}\label{eq 3}
   C(x, t) = \int_0^x \frac{\partial B}{\partial t}(\eta, t) \frac{d\eta}{(x-\eta)^\alpha}.
\end{equation}
This may be written as
\begin{equation} \label{eq 4}
C(x,t)=\Gamma(1-\alpha)J^{1-\alpha}F(x;t).
\end{equation}
Considering $t$ to be a parameter and then taking the Caputo derivative of order $1-\alpha$ of both sides of (\ref{eq 4}), we have by Eq.~(\ref{prop1})
\[
F(x;t)=\frac{1}{\Gamma(1-\alpha)}D^{1-\alpha}C(x; t).
\]
Upon integration and using the definition of the Caputo derivative we arrive at 
\begin{align}
   B(x, t) &= \int_0^t F(x, \tau)\,d\tau \nonumber \\
	         &= \frac{1}{\Gamma(1-\alpha)\Gamma(\alpha)}\int_0^t \int_0^x \frac{\partial C}{\partial x}(\xi, \tau) \frac{d\xi}{(x-\xi)^{1-\alpha}}\,d\tau  \label{caputoeq4} \\
					 &=\frac{\sin (\pi \alpha)}{\pi}\int_0^t \int_0^x \frac{\partial C}{\partial x}(\xi, \tau) \frac{d\xi}{(x-\xi)^{1-\alpha}}\,d\tau. \nonumber
\end{align}
This completes the proof.
\end{proof}
\begin{remark}
As assumed in the lemma above, $C(x,t)$ needs to be a \textit{linear combination} of separable functions of $x$ and $t$ and functions of only $x$ (but, not constant functions). So, the candidates for $C$ includes functions such as $\frac{x^{4/3}}{\sqrt{1-t^2}}$ and $\cos(x-t)+3\sqrt{x+5}$, but not functions such as $x^3+2$ or $x+t^2$. This assumption is needed in Eq.~(\ref{caputoeq4}) to guarantee that the solution obtained for $B(x,t)$ satisfies the original equation (\ref{eq2}). This is because, if there were a term that involves a function of $t$ along or constant, then that would be canceled out by taking $\frac{\partial C}{\partial x}$ and cannot be recovered by integrating with respect to $t$. 
\end{remark}
\begin{remark} This class of integral equations are normally solved by multiplying by an auxiliary function (special kernel) and then integrating (see \cite[p.30]{samko}). In our method, we used the inverse property of the Caputo derivative, that is Eq.~(\ref{prop1}), thus solving the integral equation in a simpler manner. 
\end{remark}

We also propose a new method for solving  a class of partial integro-differential equations.  This method relies on transforming the integro-differential equation into an integral equation, and then using Picard iterations to arrive at the solution.  This approach is motivated by \cite{Loverro}, where Loverro solves fractional integral equations by Picard iterations.  Here we extend his approach to partial integro-differential equations of fractional order.  

\begin{theorem}\label{thm 2}
Let $\mathcal{L}_t$ be a second order linear partial differential operator of the form
\[
\mathcal{L}_t=\frac{\p^2 }{\p t^2}+\phi(t)\frac{\p}{\p t}+\theta(t), 
\] 
for some differentiable functions $\phi,\theta$.  Let $y(t)$ satisfy
\begin{equation}\label{homg eq}
\mathcal{L}_ty=0,
\end{equation}
with $y(t)\not = 0$ for all $t$.  Then if $f(x,t)$ is continuous on $[0,1]\times [0,T]$, the problem 
\begin{equation}\label{gen prob}
\mathcal{L}_tB=J^\alpha\frac{\p B}{\p t}+f, \quad \mbox{with} \quad B(x,0)=\frac{\p B}{\p t}(x,0)=0
\end{equation}
on $[0,1]\times [0,T]$ has a unique solution given by
\begin{equation}\label{gen prob: solution}
B(x,t)=\sum_{n=0}^\infty y\left(K_2\right)^n(K _1 f ),
\end{equation}
 where the integral operators $K_1,K_2,$ and $g$ are defined by
\begin{align*}
&K_1f:=Ig(-\phi,-y)Ig(\phi,y)y^{-1}J^\alpha f,\\
&K_2K_1f:=Ig(-y,-\phi)(J^\alpha g K_1f-J^\alpha I\frac{\p(gy^{-1})}{\p t}(K_1f)),\\
&g(y,\phi):=\exp\left(I\left(2\frac{d y }{d t}y^{-1} +\phi\right)\right),
\end{align*} 
and $(K_2)^n$ denotes the operator composition $K_2^nK_1f:=K_2(K_2(\cdots (K_1f)\cdots))$, $n$ times.

\end{theorem}
\begin{proof}
 To prove existence we will transform equation (\ref{gen prob}) into a linear integral equation, and then use an inductive (iterative) process to arrive at the solution. In order to do this, an integrating factor method would be helpful. However, with the equation in its current form, we are unable to do so due to the presence of the term $\theta B$.  To eliminate this we will proceed with an {\it ansatz} of the form
\[
B(x,t)=b(x,t)y(t),
\]
where $y$ solves the homogenous equation (\ref{homg eq}) and $y(t)\not =0$ for any  $t \ne 0$.  Substituting this into (\ref{gen prob}) we obtain
\begin{equation*}
\left( y\frac{\p^2 b}{\p t^2}+2\frac{\p b}{\p t}\frac{d y }{d t}+b\frac{d^2 y}{d t^2} \right)+\phi\left( y\frac{\p b}{\p t} +b\frac{d y}{d t}\right)+\theta by=J^\alpha \frac{\p }{\p t}(yb)+f.
\end{equation*}
However since $y$ is a solution of $\mathcal{L}_t=0$, we have
\begin{equation*}
 \frac{\p^2 b}{\p t^2}+\left( 2\frac{d y }{d t}y^{-1}+\phi \right)\frac{\p b}{\p t} =y^{-1}\left(J^\alpha \frac{\p }{\p t}(yb)+f\right).
\end{equation*}
This is a Bernoulli-type linear differential equation for the term $\frac{\p b}{\p t}$. Thus, following the standard procedure and multiplying by an integrating factor we have
\begin{equation}\label{pre int factor}
\frac{\p}{\p t}\left(\exp\left(I\left( 2\frac{d y }{d t}y^{-1} +\phi\right) \right)\frac{\p b}{\p t} \right)= \exp\left(I\left( 2\frac{d y }{d t}y^{-1} +\phi\right)\right)y^{-1}\left(J^\alpha \frac{\p }{\p t}(yb)+f\right).
\end{equation}
Before proceeding, for simplicity, let
\begin{align*}
&g(y,\phi):=\exp\left(I\left(2\frac{d y }{d t}y^{-1} +\phi\right)\right)\\
\end{align*}
and
\begin{align*}
&K_1f:=Ig(-\phi,-y)Ig(\phi,y)y^{-1}f,\\
&K_1J^\alpha\frac{\p}{\p t}(yb):=Ig(-\phi,-y)Ig(\phi,y)y^{-1}J^\alpha\frac{\p}{\p t}(yb).
\end{align*}
Then integrating (applying $I$ to) each side of (\ref{pre int factor}), multiplying each side by $g(-\phi,-y)$, and integrating again we arrive at
\begin{equation*}
b=K_1\left( J^\alpha\frac{\p}{\p t}(yb)+f\right).
\end{equation*}
In order to get the time derivative off of the term $\frac{\p (yb)}{\p t}$, we may apply integration by parts after applying Fubini's theorem. 
 Indeed, since $gy^{-1}$ is independent of $x$:
\[
Igy^{-1}J^\alpha=IJ^\alpha gy^{-1}.
\]
Therefore if $\frac{\p B}{\p t}=\frac{\p (by)}{\p t}$ is integrable (which we will show later), 
\[
\frac{1}{(x-\nu)^{1-\alpha}}{\frac{\p B}{\p t}}
\]
is integrable on $[0,T]\times [0,1]$.  So we may apply Fubini's theorem
\[
IJ^\alpha gy^{-1}\frac{\p (yb)}{\p t}=J^\alpha Igy^{-1}\frac{\p (yb)}{\p t},
\]
and  integration by parts to obtain
\[
J^\alpha Igy^{-1}\frac{\p (yb)}{\p t}=J^\alpha g b-J^\alpha I\frac{\p(gy^{-1})}{\p t}(yb).
\]
Now we define $K_2b$ by
\begin{align*}
K_2b:=Ig(-y,-\phi)\left(J^\alpha g b-J^\alpha I\frac{\p(gy^{-1})}{\p t}(yb)\right),
\end{align*}
in order to cast our differential equation as an integral equation:
\begin{equation}\label{integral equation}
b=K_2b+K_1f.
\end{equation}
Similarly we define the iterative sequence
\[
b_{n+1}=K_2b_n+K_1f.
\]
Then for some initial guess $b_0$  we have
\begin{equation}\label{it_seq}
b_{n+1}=\sum_{i=0}^{n}(K_2)^iK_1f+(K_2)^{n+1}b_0.
\end{equation}
If the sequence defined by (\ref{it_seq}) converges (see the Appendix), then we have
\begin{equation}\label{b series}
b=\sum_{i=0}^{\infty}(K_2)^iK_1f,
\end{equation}
which gives rise to (\ref{gen prob: solution}) upon substitution of $B=yb$ into the product.  Also observe that (\ref{integral equation}) implies uniqueness, for the associated homogenous problem is a linear integral equation for which 0 is a solution.
\end{proof}

\begin{remark}
The assumption $B(x,0)=0$ is not necessary, and was just made to simplify the algebra.  The theorem applies with more general initial conditions such 
as $B(x,0)=h(x),\frac{\p B}{\p t}(x,0)=0$; in those cases the form of (\ref{gen prob: solution}) would simply change.  Also note the 
second-order differential operator $\mathcal{L}_t$ need not have constant coefficients.  We require only the knowledge of a solution to the corresponding equation $\mathcal{L}_ty=0$,  
with $y(t)\not = 0$. 

\end{remark}

We have proved the theorem when $\mathcal{L}_t$ is a linear non-constant-coefficient, second-order differential operator.  But similar results hold when $\mathcal{L}_t$ is first order, or even when $\mathcal{L}_t$ is a temporal Riemann--Liouville fractional derivative.  Indeed,
if we have
\[
D^\beta_t B=J^\alpha B +f,
\]
then 
\[
B=J^\beta_tJ^\alpha B+J^\beta_tf,
\]
which is a linear integral equation, and hence may be solved by Picard's iteration method.  The solution is given as
\[
B=\sum_{n=0}^\infty( J^\beta_tJ^\alpha )^n J^\beta_t f.
\]
This may be summarized in the following theorem.

\begin{theorem}\label{thm 3}
Let $B,f:[0,1]\times [0,T]\to \mathbb{R}$ be differentiable with respect to $t$, and integrable with respect to $x$. 
Then the problem  
\begin{equation}\label{frac diff}
D^\beta_t B=J^\alpha B +f, \qquad B(x,0)=f(x,0)=0
\end{equation}
has a unique solution given by
\[
B=\sum_{n=0}^\infty( J^\beta_tJ^\alpha )^nJ^\beta_t f.
\]
\end{theorem}
\begin{proof}
We first apply the fractional integral operator $J_t^\beta$ to each side of (\ref{frac diff}).  This transforms (\ref{frac diff}) into
\begin{equation}
B=J_t^\beta J^\alpha B+J_t^\beta f.\label{rl temp}
\end{equation}
In a manner analogous to the proof of Theorem \ref{thm 2}, we define the iterative sequence  
\begin{equation}
B_{n+1}=J_t^\beta J^\alpha B_n+J_t^\beta f,\label{temp it}
\end{equation}
so that the solution to (\ref{rl temp}) is given by $B(x,t)=\lim_{n\to\infty}B_n(x,t)$.  Then given some initial function $B_0$,  the sequence (\ref{temp it}) implies 
\[
B_{n+1}=\sum_{i=0}^n (J_t^\beta J^\alpha)^i(J_t^\beta f)+(J_t^\beta J^\alpha)^{n+1}B_0.
\]
Hence, the solution to (\ref{rl temp}) is given by
\[
B=\lim_{n\to\infty}B_n=\sum_{i=0}^\infty (J_t^\beta J^\alpha)^i(J_t^\beta f).
\]
Convergence may be verified by an application of Weierstrass's $M$-test, in a manner similar to the one given in the Appendix.
\end{proof}

\section{Illustrative Examples}\label{Section: Illustrative Examples}

\begin{example}As the first example, consider the following Abel-type multivariate integral equation (of the first kind) with a homogeneous initial condition: 
\begin{equation}\label{ex1}
    \frac{x^{4/3}}{\sqrt{1-t^2}}= \int_0^x \frac{\partial B}{\partial t}(x-\xi, t) \frac{d\xi}{\xi^{\frac{2}{3}}} \quad \mbox{with} \quad B(x, 0) = 0.
\end{equation}
According to Lemma~\ref{lemma1}, taking $C(x, t) = \frac{x^{4/3}}{\sqrt{1-t^2}}$, it is straightforward to see that the solution of (\ref{ex1}) is given by 
\begin{align}
B(x, t) &= \frac{\sin (\frac{2\pi}{3} )}{\pi}\int_0^t\int_0^x \frac{\partial C}{\partial x}(\xi, \tau) \frac{d\xi}{(x -\xi)^{\frac{1}{3}}}\,d\tau \nonumber \\
        &= \frac{\sin (\frac{2\pi}{3} )}{\pi}\int_0^t\int_0^x \frac{4}{3}\xi^{\frac{1}{3}}\frac{d\xi}{(x -\xi)^{\frac{1}{3}}}\frac{d\tau}{\sqrt{1-\tau^2}} \nonumber \\
				&= \frac{4x}{3}\cdot\frac{\sin (\frac{2\pi}{3} )}{\pi}\int_0^t\frac{d\tau}{\sqrt{1-\tau^2}}\int_0^1 u^{\frac{1}{3}}(1 -u)^{-\frac{1}{3}}\,du, \quad \mbox{where} \quad u = \frac{\xi}{x},\nonumber \\
				&= \frac{4}{9} x\sin^{-1}(t). \label{e3}
\end{align} Notice that we have used the properties of the $\Gamma(\cdot)$ function such as $$
B(\alpha, \beta) = \frac{\Gamma(\alpha)\Gamma(\beta)}{\Gamma(\alpha+\beta)}=\int^1_0 u^{\alpha -1}(1-u)^{\beta-1} du,\quad  \Gamma(n+1) = n\Gamma(n),\quad \Gamma(\alpha)\Gamma(1-\alpha)=\frac{\pi}{\sin (\pi\alpha)}.
$$ 
\end{example}
Another advantage of the method we proposed in this paper is that it can even be applied to certain Abel-type integral equations with inhomogeneous initial conditions. For example, consider an initial condition of the form $B(x, 0) = f(x)$, where $f(x)$ is a differentiable function. In this case, we can define a new function, $D(x, t) = B(x, t) - f(x)$. Then, we see that $\frac{\partial D}{\partial t}=\frac{\partial B}{\partial t}$ and $D(x, 0)=0$. Thus, we can first, solve the integral equation for $D(x, t)$ using Lemma~\ref{lemma1} and then find $B(x,t) = D(x, t) + f(x)$. To see this, consider the following problem instead.
\begin{equation}\label{ex2}
    \frac{2}{3}\cos(x-t)= \int_0^x \frac{\partial B}{\partial t}(x-\xi, t) \frac{d\xi}{\xi^{\frac{2}{3}}} \quad \mbox{with} \quad B(x, 0) = e^{\tan(x)}+4.
\end{equation}
First let $D(x,t)=B(x,t)-e^{\tan(x)}-4$. Now, using Lemma~\ref{lemma1} and taking $C(x, t) = \frac{2}{3}\cos(x-t)$, we have
\begin{align}
D(x, t) &= \frac{\sin (\frac{2\pi}{3} )}{\pi}\int_0^t\int_0^x \frac{\partial C}{\partial x}(\xi, \tau) \frac{d\xi}{(x -\xi)^{\frac{1}{3}}}\,d\tau \nonumber \\
        &= -\frac{2\sin (\frac{2\pi}{3} )}{3\pi}\int_0^t\int_0^x \sin(\xi-\tau)\frac{d\xi}{(x -\xi)^{\frac{1}{3}}}\,d\tau \label{neqn1} \\
				&= \frac{\sqrt{3}}{\pi}x^{2/3}(1-\cos(t)){}_1F_2\left(1;\frac{5}{6},\frac{4}{3};-\frac{x^2}{4}\right)     
				      -\frac{3\sqrt{3}}{10\pi}x^{5/3}\sin(t){}_1F_2\left(1;\frac{4}{3},\frac{11}{6};-\frac{x^2}{4}\right), \quad x>0, \label{neqn2}
\end{align}
where ${}_1F_2(a;b,c;x)$ is the Gauss Hypergeometric function. The result in Eq.~(\ref{neqn2}) has been obtained using the CAS Mathematica\textsuperscript{\textregistered}.  
Thus, the solution of the initial value problem is
\[
   B(x,t) = \frac{\sqrt{3}}{\pi}x^{2/3}(1-\cos(t)){}_1F_2\left(1;\frac{5}{6},\frac{4}{3};-\frac{x^2}{4}\right)     
				      -\frac{3\sqrt{3}}{10\pi}x^{5/3}\sin(t){}_1F_2\left(1;\frac{4}{3},\frac{11}{6};-\frac{x^2}{4}\right) + e^{\tan(x)}+4, \quad x>0. 
\]
It can be seen from this example that even a simple looking integral equation may be quite involved sometimes. The integral in Eq.~(\ref{neqn1}) can be also evaluated without the hypergeometric functions approach as follows. To that end, notice that
\begin{align}
     \int_0^t\int_0^x \frac{\sin(\xi-\tau)}{(x -\xi)^{\frac{1}{3}}}\,d\xi\,d\tau &= \int_0^t\left(\cos \tau \int_0^x\frac{\sin \xi}{(x-\xi)^{\frac{1}{3}}}\,d\xi -\sin \tau \int_0^x\frac{\cos \xi}{(x-\xi)^{\frac{1}{3}}}\,d\xi \right)d\tau  \nonumber \\
    &=\sin t\int_0^x\frac{\sin \xi}{(x-\xi)^{\frac{1}{3}}}\,d\xi+(\cos t -1)\int_0^x\frac{\cos \xi}{(x-\xi)^{\frac{1}{3}}}\,d\xi \nonumber \\
		&=\frac{3}{2}\sin t\int_0^x (x-\xi)^{\frac{2}{3}}\cos \xi\,d\xi+(\cos t -1)\left[\frac{3}{2}x^{\frac{2}{3}}-\frac{3}{2}\int_0^x(x-\xi)^{\frac{2}{3}}\sin \xi\,d\xi\right] \nonumber \\	
		&= \frac{3}{2}\sin t\int_0^x u^{\frac{2}{3}}\cos (x-u)\,du+(\cos t -1)\left[\frac{3}{2}x^{\frac{2}{3}}-\frac{3}{2}\int_0^x u^{\frac{2}{3}}\sin (x-u)\,du  \right] \nonumber \\	
		&= \frac{3}{2}(\cos t -1)x^{\frac{2}{3}}+\frac{3}{2}\left[\sin(t-x)+\sin x \right]\int_0^x u^{\frac{2}{3}}\cos u\,du \nonumber \\
		&\hspace{3.5cm} + \frac{3}{2}\left[\cos(t-x)-\cos x \right]\int_0^x u^{\frac{2}{3}}\sin u\,du  \nonumber			
\end{align}
Now, using Taylor expansions for $\cos u$ and $\sin u$, we can see that 
\[
    \int_0^x u^{\frac{2}{3}}\cos u\,du = x^{\frac{5}{3}}\sum_{k=0}^\infty (-1)^k\frac{x^{2k}}{(2k+\frac{5}{3})(2k)!},\quad\mbox{and} \quad\int_0^x u^{\frac{2}{3}}\sin u\,du = x^{\frac{8}{3}}\sum_{k=0}^\infty (-1)^k\frac{x^{2k}}{(2k+\frac{8}{3})(2k+1)!}.
\]
Thus, the solution of the integral equation in Eq.~\ref{ex2}, takes the form
\begin{align}
   B(x,t) = -\frac{\sqrt{3}}{2\pi}\Big[
	(\cos t -1)x^{\frac{2}{3}} &+ \left(\sin(t-x)+\sin x \right)x^{\frac{5}{3}}\sum_{k=0}^\infty (-1)^k\frac{x^{2k}}{(2k+\frac{5}{3})(2k)!} \nonumber \\
	&+ \left(\cos(t-x)-\cos x \right)x^{\frac{8}{3}}\sum_{k=0}^\infty (-1)^k\frac{x^{2k}}{(2k+\frac{8}{3})(2k+1)!}\Big] + e^{\tan(x)}+4, \quad x>0. \nonumber
\end{align}
Notice that, in the proof, we used the substitution $u=x-\xi$, integration-by-parts method, some trigonometric identities, and Taylor expansions. 

Next we will study some applications of our method that arise in modeling surface-volume reactions. 

\section{Applications} \label{Applications}

\subsection{Setting}\label{Setting}

One of the major goals of this paper is to explore the applications of fractional calculus in solving real-world applications that arise in science, engineering and other disciplines. For example, stereology is concerned with the determination of three-dimensional objects from two- or one-dimensional data, and the extrapolation to three dimensions is based on the solution of a mathematical model derived using geometric probability and statistics.  It is quite often the case that the model emerges in the form of an Abel integral equation \cite[p.~123]{jake} given by
\[
\int_{\rm a}^\infty \left(\frac{x}{a}\right)^{\mu -1} \frac{g(\sqrt{x/a})}{(x-a)^\mu}\, dx = kz(a),
\]
where $g(x)$ and $z(x)$ are probability density functions. 

The Abel integral equations of the first/second kind with $\alpha = \frac{1}{2}$ also appear in a variety of problems in physics and chemistry; see, for example, the references in \cite{brun2,olm,jake}. 
In this section we discuss a situation where, for the first time, we come across a fractional order different from $\frac{1}{2}$ in a real-world application. 

The applications we consider here arise in the setting of \emph{surface-volume reactions}.  A surface-volume reaction is one where a buffer fluid containing ligand molecules is convected through a channel over a surface to which immobilized ligands (or receptors) are confined 
(see Figure \ref{cross_sect}).  The scope of surface-volume reactions is very broad, and are of great importance in biology.  Such reactions occur in blood clotting \cite{Grabowski}, drug absorption \cite{Bertucci}, and antigen-antibody  interactions \cite{Raghaven}.  Surface-volume reactions also occur in DNA-protein interaction, which affects gene expression \cite{Edwards1}.  Additionally, purification processes often occur in channels with reactants embedded along the wall  \cite{Edwards1}.  In order to experimentally study such reactions, scientists use an apparatus known as an \emph{optical biosensor}.  The use of optical biosensors is quite popular, with over 10,000 authors citing the use of an optical biosensors as of 2009 alone \cite{Rich}.  A schematic is depicted below.
\begin{figure}[H]
\centering
  \includegraphics[width=13.7cm, height=7.5cm]{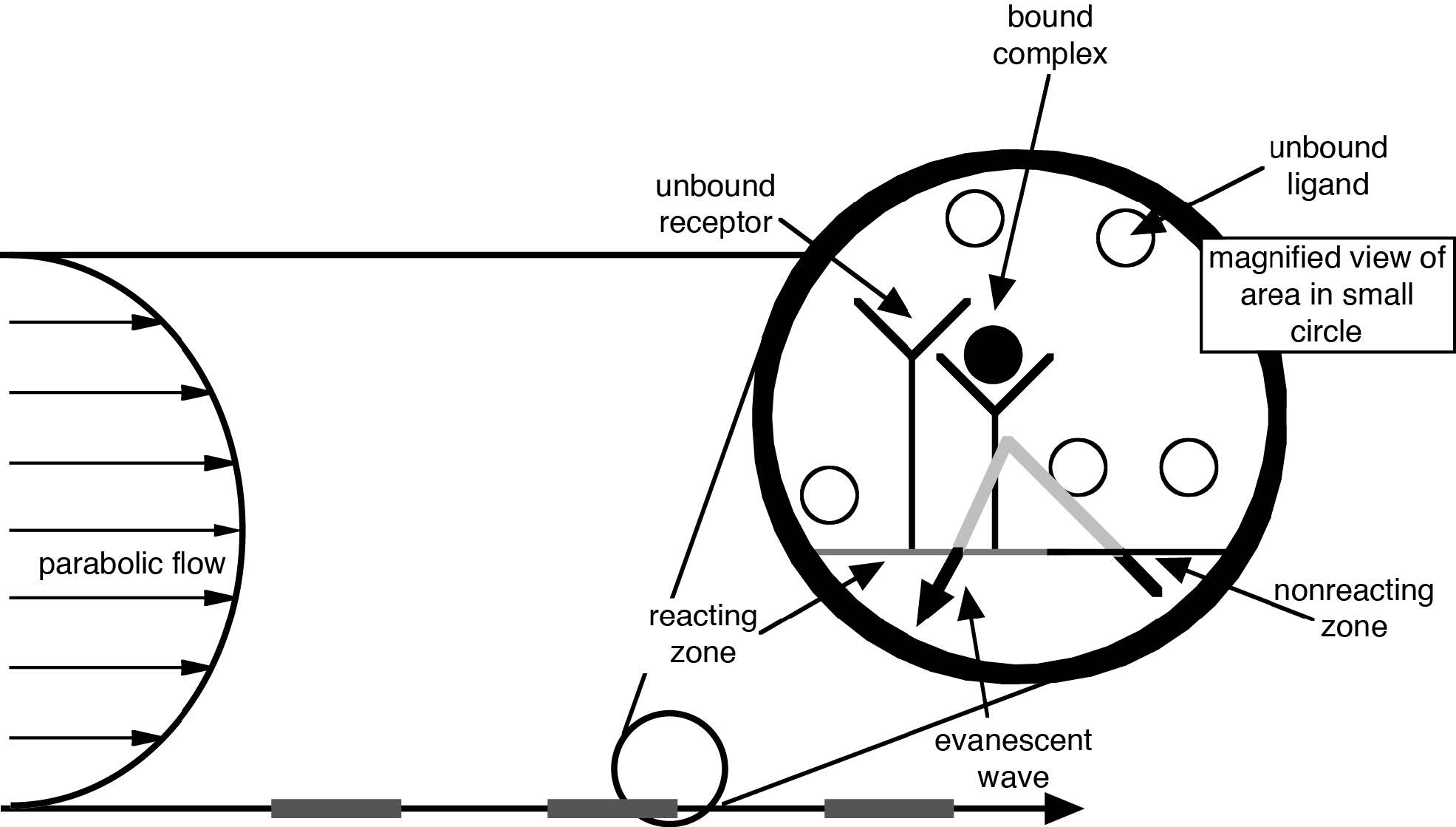}
  \caption{Cross-sectional schematic of the channel and binding/unbinding. Taken from \cite{Edwards2}.}
	\label{cross_sect}
\end{figure}
In Figure \ref{cross_sect} one can see the unbound ligand being convected through the channel.  As the unbound ligand molecules bind with the receptors on the surface, an evanescent wave is reflected off of the floor of the device.  Refractive changes due to binding of the reactants are then averaged over the length of the ceiling to provide scientists 
with real-time mass measurements of the bound ligand concentration.   The chemical kinetics at the boundary may be written as 

\[E+L \xrightleftharpoons[\wt{k}_{\rm d}]{\wt{k}_{\rm a}}L,\] 
where $E$ represents an empty receptor and $L_i$ denotes an unbound ligand molecule at the surface.  A single ligand may bind with each receptor.  The parameter $\wt{k}_{\rm a}$ represents the dimensional reaction rate of ligand binding, and $\wt{k}_{\rm d}$ represents the dimensional dissociation rate.  Herein we will denote dimensional quantities with a tilde.  Also we denote the dimensionless concentration of bound receptor sites at the surface by
\begin{align*}
&[EL]=B,\ [L]=C.
\end{align*}
In general $C=C(x,y,t)$, and because there are only receptor sites on the floor of the channel, $B$ is a function of $x$ and $t$ only. Following \cite{Edwards1} we have scaled both the dimensional bound ligand concentration $\tilde{B}$ and the dimensional unbound concentration $\wt{C}$, so that the dimensionless variables $B,C$ are between 0 and 1.  Since there are no other sources or sinks, the bound ligand will change only due to association or dissociation:
\begin{equation}
\frac{\p B}{\p t}=(1-B)C(x,0,t)-KB, \qquad B(x,0)=0.
\label{kinetics}
\end{equation}
Here $K$ is a dimensionless parameter containing the ratio of the association and dissociation coefficients $\wt{k}_{\rm a}, \wt{k}_{\rm d}$.  The term $(1-B)C$ represents a bimolecular production term, and $-KB$ represents dissociation. If the unbound ligand concentration at the boundary $C(x,0,t)$ were uniform, taking $C=1$ gives the following solution of (\ref{kinetics}):
\begin{align*}
&B(t)=\rho^{-1}(1-e^{-\rho t}),\qquad \rho=K+1.
\end{align*}
There is initially no bound ligand in the channel, and as time progresses ligand molecules bind with the receptors until a balance between association and dissociation is reached, and the system reaches chemical equilibrium.  Due to dissociation effects, this will happen before saturation.

\begin{figure}
\centering
  \includegraphics[width=9.0cm, height=5.5cm]{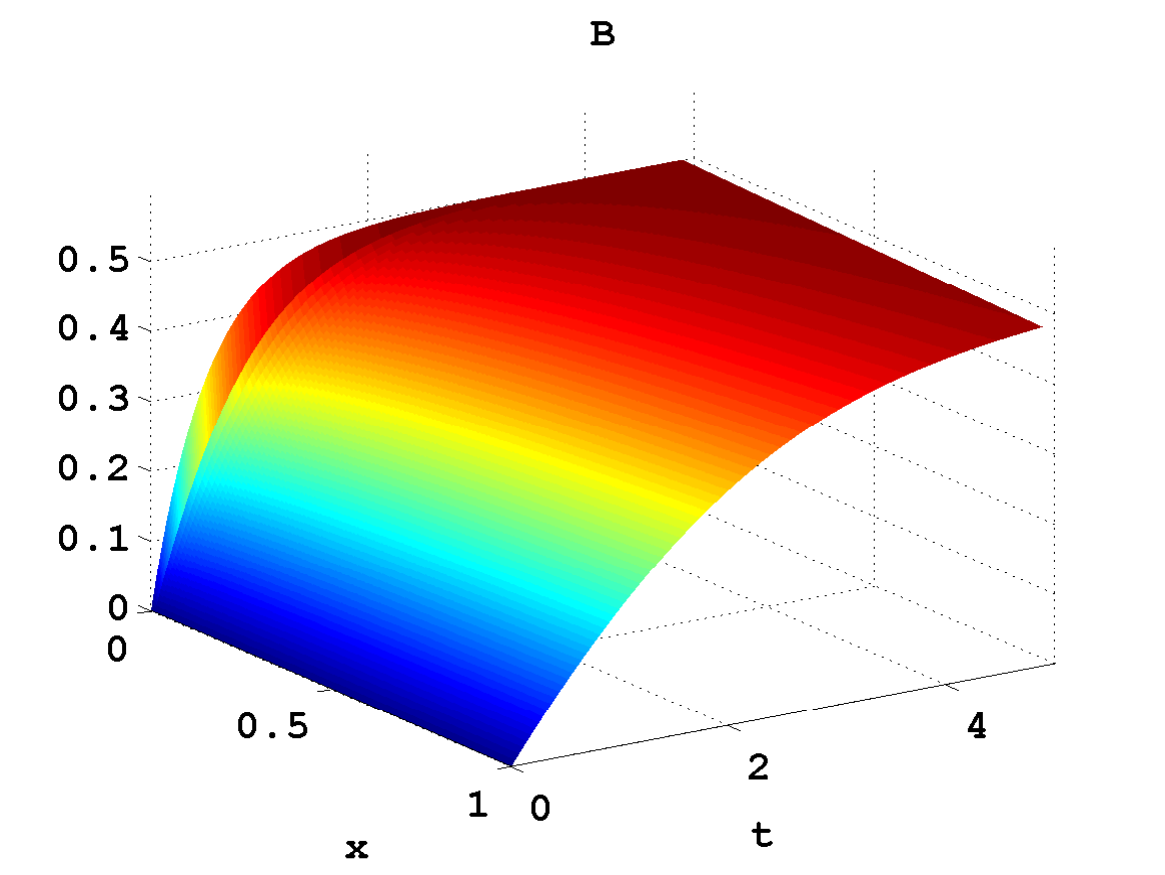}
  \caption{Numerical solution of (\ref{main eq}), $\Da=2$, $K=1$, computed numerically using the algorithm described in \cite{Edwards3}.}
    \label{fig: num soln}
\end{figure}

In reality, transport will not be perfectly efficient.  Edwards \cite{Edwards1} has shown that $C(x,0,t)$ takes the form of a fractional integral
\begin{align*}
C(x,0,t)&=1-\frac{\Da}{3^{1/3}\Gamma(2/3)}\int_0^x\frac{\p B}{\p t}(\nu,t)(x-\nu)^{-2/3}d \nu,
\end{align*}
which may be expressed in terms of the Riemann-Louiville fractional integral
\begin{align}
C(x,0,t)=&1-\frac{\Da\Gamma(1/3)}{3^{1/3}\Gamma(2/3)}J^{1/3}\frac{\p B}{\p t}.\label{frac int}
\end{align}
Substituting $C(x,0,t)$ into the kinetics equation we arrive at
\begin{equation}
\frac{\p B}{\p t}=(1-B)\left(1-\frac{\Da\Gamma(1/3)}{3^{1/3}\Gamma(2/3)}J^{1/3}\frac{\p B}{\p t}\right) -KB, \qquad B(x,0)=0,
 \label{main eq}
\end{equation}
where $\Da$ denotes the Damk{\"o}hler number, an important dimensionless parameter representing a ratio of reaction to diffusion. The Damk{\"o}hler number
\[
\Da=\frac{\wt{k}_{\rm a} \wt{L}^{1/3}\wt{h}^{1/3}\wt{R}_{\rm t}}{\wt{V}^{1/3}\wt{D}^{2/3}}
\]
 is a function of the dimensions of the channel $\wt{h}, \wt{L}$,  characteristic velocity scale $\wt{V}$, diffusion rate of the ligand molecules $\wt{D}$,  the reaction rate $\wt{k}_{\rm a}$, and the concentration of total receptor sites $\wt{R}_{\rm t}$.  For physically realizable scenarios, $\Da$ is either $o(1)$ or $O(1)$ \cite{Edwards1}.
  In Figure \ref{fig: num soln} we took $\Da=2$ to exaggerate transport effects; however in our analysis we will consider only the experimentally relevant regime, $\Da=o(1)$. 
 This corresponds to the parameter regime in which diffusion is much quicker than reaction; \textit{i.e.}, when transport is quite efficient. Note that (\ref{kinetics}) 
is recovered in the limit as the Damk{\"o}hler number goes to zero in (\ref{main eq}). For a visualization of transport effects on the evolution of the bound ligand concentration, 
see Figure \ref{fig: num soln}.

\subsection{Physical Interpretation} 

We now discuss the physical interpretation of the fractional integral in (\ref{frac int}).  There are multiple time scales associated with the problem.  There is a time scale for convection $t_{\rm c}$,  diffusion near the wall $t_{\rm w}$, diffusion into the surface $t_{\rm d}$, and reaction $t$.   The flow \emph{away from the wall} reaches equilibrium on the $t_{\rm c}$ time scale, and reaction occurs on the latter of the two time scales.  The time scale for reaction is much slower than convection.  Thus, one would expect that the unbound concentration at the boundary $C(x,0,t)$ would be a perturbation away from the uniform outer concentration of 1:
\begin{equation}
C(x,0,t)=1-\Da c(x,t).
\label{pert exp}
\end{equation}
Edwards showed in \cite{Edwards1} that the perturbation $c(x,t)$ is given by
\[
c(x,t)=\frac{\Gamma(1/3)}{3^{1/3}\Gamma(2/3)}J^{1/3}\frac{\p B}{\p t}
,
\]
which gives (\ref{frac int}) upon substitution into (\ref{pert exp}).

When transport effects are considered, initially there will be slightly more bound ligand upstream than downstream.  
This is due to the fact that the ligand will diffuse into the surface upstream first.  Thus the fractional integral in (\ref{frac int}) represents 
\emph{upstream ligand depletion} \cite{Edwards1}.  If $x$ is smaller, then less ligand will have  already bound, thus (\ref{frac int}) will be larger, and
 there will be more ligand available for binding at the surface upstream.  Note that $B$ is increasing so $\frac{\p B}{\p t}$ will be non-negative.

Thus far we have discussed only the \emph{association phase} (or the \emph{injection phase}) of the experiment.  One can also model the \emph{dissociation} 
or\emph{ wash phase} of the experiment.  After studying the association phase of the experiment, scientists may wish to clean the optical biosensor for reuse.
  Therefore to clean the device scientists will simply convect only the buffer fluid through the channel.  This has the effect of washing out all of the bound ligand. 
 In the absence of transport this can be well modeled by the standard exponential decay curve
\[
 B(t)=\rho^{-1}e^{-\rho t}.
 \]
This follows from (\ref{kinetics}) with $C=0$ and $B(0)=\rho^{-1}$ (the equilibrium concentration in the association phase).  However when considering transport effects, Edwards has shown in \cite{Edwards2} that the kinetics process is governed by
\[
\frac{\p B}{\p t}=(1-B)\left[\frac{-\Da}{3^{1/3}\Gamma(2/3)}\int_0^x\frac{\p B}{\p t}(x-\nu)^{-2/3}\,d \nu\right] -KB, \qquad B(x,0)=\rho^{-1}, 
\] 
or
\begin{equation}
\frac{\p B}{\p t}=(1-B)\left(-\frac{\Da\Gamma(1/3)}{3^{1/3}\Gamma(2/3)}J^{1/3}\frac{\p B}{\p t}\right) -KB.
\label{diss}
\end{equation}
Note that in this phase $B$ is decreasing, so $\frac{\p B}{\p t}$ will be negative and the term
\[
(1-B)\left[\frac{-\Da}{3^{1/3}\Gamma(2/3)}\int_0^x\frac{\p B}{\p t}(x-\nu)^{-2/3}\,d \nu\right]
\]
will be positive.  Observe that if $\Da=0$, transport is perfectly efficient and we recover the standard exponential decay model.  The fractional integral in this phase has a slightly 
different interpretation.  Ligand will  bind in the wash phase only when ligand molecules dissociating upstream rebind with receptor sites downstream.  If a receptor site is further 
downstream, then it is more likely that a ligand that has dissociated upstream will rebind with it.  This is exactly what the fractional integral in (\ref{diss}) tells us.  
When $x$ is very small, the fractional integral in (\ref{diss}) will also be small, and there will be a small possibility that a ligand molecule will rebind. 
 Near the end of the channel, $x$ will be close to 1, and inefficient transport will result in rebinding.  Thus in this phase the fractional integral accounts for 
ligand rebinding due to inefficient transport.  Notice that when $\Da$ is larger, the fractional integral will have a larger effect, since transport is less effecient. 
 If transport is less efficient than the ligand rebinding effect will be more exaggerated, which is what (\ref{diss}) tells us.

\subsection{Applications of Theorem \ref{thm 2}}

\subsubsection{A Multiple Scale Expansion}

In the experimentally relevant parameter regime of $\Da\ll1$ one may propose a perturbation series of the form
\begin{equation}
B(x,t)=B_0(t)+\Da B_{1}(x,t)+O(\Da^2),\label{B pert}
\end{equation}
to investigate the behavior of the solution of (\ref{main eq}).  Here the $B_i$ are assumed to be independent of $\Da$.  Such a regular expansion may be shown to be secular, motivating a search for a multiple-scale solution of the form:
\begin{align*}
&B(x,t)=B_0(x,T,\tau)+\Da B_{1}(x,T,\tau)+O(\Da^2),
\end{align*}
where
\[ T=\Da t \quad \mbox{and} \quad \tau=\left(1+\sum_{n=2}^\infty \omega_n\Da^n\right) t. \]
Here the $\omega_n$ are constants used to suppress higher-order secularities.  

When doing so, in order to eliminate a secular term, one must solve the equation \cite{Edwards1}
\begin{equation}\label{eqn:feqn}
\frac{\p a}{\p t}=\frac{K}{3^{1/3}\Gamma(2/3)}\int_0^x a(\nu,t)(x-\nu)^{-2/3}\, d\nu,\qquad \ a(x,0)=-\frac{1}{\alpha},\qquad x\in[0,1].
\end{equation}

In \cite{Edwards1}, Edwards extends the domain of \eqref{eqn:feqn} to $x\in[0,\infty)$, relying upon physically realistic arguments.  Once that assumption has been made, \eqref{eqn:feqn} can be solved using Laplace transforms in $x$.  The solution in transform space is expanded in a Taylor series, and then inverted term by term to find $a$.  The solution thus obtained is used to obtain $B(x,t)$
on the domain of physical interest by restriction to $x\in[0,1]$.

Though perhaps justifiable physically, such an approach is mathematically unsatisfying, as it relies on an extension of the domain.  We can also use fractional calculus methods to get the same result on the original domain $x\in[0,1]$.  First let us rewrite (\ref{eqn:feqn}) using the Riemann--Liouville integral:
\begin{align*}
&\frac{\p  a}{\p t}=rJ^{1/3}a, \quad \mbox{where} \; r=\frac{\Gamma(1/3)K}{\Gamma(2/3)3^{1/3}}.
\end{align*}
Integrating both sides gives the integral equation
\begin{equation}
 a =r\int_0^t J^{1/3}a\, d t-\rho^{-1}.\label{ms integral equation}
\end{equation}
By analogy with the proofs of Theorems \ref{thm 2} and \ref{thm 3} we can define the iterative sequence
\begin{equation*}
a_{n+1}=rIJ^{1/3}a_n-\rho^{-1}, \quad a_0=-\rho^{-1}.
\end{equation*}
Hence, in a manner similar to the proof of Theorems \ref{thm 2} and \ref{thm 3}, one may show $a_n$ converges to the solution of (\ref{ms integral equation}):
\begin{equation}
\label{ eq 2}
 a(x,t)=-\rho^{-1}\sum_{n=0}^\infty\frac{(r tx^{1/3})^n}{\Gamma(1+n/3)n!}. 
\end{equation}
The above is precisely the solution Edwards obtains in \cite{Edwards1}. It is interesting to see that the expression in \eqref{ eq 2} is nothing but 
the Hadamard product of the two functions $e^{rt}$ and the Mittag-Leffler function \cite{pod}
 $$
E_\frac{1}{3} (x^{1/3}) = \sum_{n = 0}^\infty \frac{(x^{1/3})^n}{\Gamma(1+n/3)}.
$$

\subsubsection{A Linearized Equation}

Equation (\ref{main eq}) is unwieldy and hopeless to solve in closed form.  However, it is of mathematical interest to consider a linear variant given by:
\begin{equation}
\frac{\p B}{\p t}+\rho B=1-\beta J^{1/3}\frac{\p B}{\p t},\qquad B(x,0)=0,\qquad \beta=\frac{\Da \Gamma(1/3)}{3^{1/3}\Gamma(2/3)}.\label{lin}
\end{equation}
The above equation does not fall under the class considered in Theorem \ref{thm 2}, only due to the fact that (\ref{lin}) is first-order in time, while the class of equations
considered in Theorem \ref{thm 2} are second-order in time.  However, we may apply the same method to solve (\ref{lin}), and begin by multiplying each side of (\ref{lin}) by the
integrating factor $e^{\alpha t}$:
\begin{align*}
&e^{\rho t}\left(\frac{\p B}{\p t}+\rho B\right)=e^{\rho t}\left(1-\beta J^{1/3}\frac{\p B}{\p t}\right),
\end{align*}
which gives
\begin{align*}
&B= \rho^{-1}(1-e^{-\rho t})-\beta I J^{1/3}e^{\rho(s-t)}\frac{\p B}{\p t},
\end{align*}
where we have used the spatial independence of the exponential function.  Next we interchange the order of integration by appealing to Fubini's Theorem:
\begin{align*}
&B= \rho^{-1}(1-e^{-\alpha t})-\beta  J^{1/3}Ie^{\rho(s-t)}\frac{\p B}{\p t}.
\end{align*}
Note an integration by parts gives
\begin{align*}
I\left(e^{\rho(s-t)}\frac{\p B}{\p t}\right)=B-\rho I e^{\rho (s-t)}B.
\end{align*}
This implies  
\begin{align*}
&B= \rho^{-1}(1-e^{-\rho t})-\beta  J^{1/3}(B-\rho I e^{\rho (s-t)}B),
\end{align*}
or rearranging,
\begin{align*}
 &[1-\beta J^{1/3}(-1 +\rho I e^{\rho(s-t)})]B= \rho^{-1}(1-e^{-\rho t}).
\end{align*}
By inverting the operator on the left-hand side, after some algebra we arrive at
\begin{align*}
B= \sum_{n=0}^\infty \frac{(\beta x^{1/3})^n}{\Gamma(1+n/3)}\left(\sum_{j=0}^n{n \choose j }(-1)^{n-j}\rho^{j}(Ie^{\rho(s-t)})^j\rho^{-1}(1-e^{-\rho s})\right).\\
\end{align*}
Here $(Ie^{\rho(s-t)})^j$ denotes the iterated $j$-fold integral
\begin{align*}
(Ie^{\rho(s-t)})^j\alpha^{-1}(1-e^{-\rho s})&=\underbrace{\int_0^t e^{\rho(s_j-t)}\cdots \int_0^{s_2}e^{\rho(s_1-s_2)}\rho^{-1}(1-e^{-\rho s_1})\ d s_1 \ldots d s_j.}_{j\ \text{times}}\\
\end{align*}
Convergence and differentiability follow from applications of the $M$-test.



\section{Conclusion}

Recently fractional calculus has been an area that is rich in application.  We have found yet another application of fractional calculus in modeling surface-volume reactions.  
In particular, we have considered  applications that occur when scientists simulate surface-volume reactions experimentally in an optical biosensor.  We have seen that when modeling
 surface-volume reactions in optical biosensors, the fractional integral may be interpreted as a \emph{ligand depletion term} during the injection phase of the experiment.  In the
 dissociation phase the fractional integral represents ligand rebinding due to inefficient transport.  Thus we have another interpretation to a centuries-old problem:
 ``what is meant by a fractional integral operator?''  This is also the first time we identify a fractional integral with order other than $\frac{1}{2}$ in an engineering problem. 

Motivated by equation (\ref{main eq}), we have considered several related equations.  We have shown that these equations may be solved using fractional differentiation and integration. 
 Additionally we have extended the approach in \cite{Loverro} to non-constant-coefficient partial integro-differential equations.  By considering equations motivated by the chemical kinetics in the surface-volume
 geometry, we have been able to extend the theory of fractional calculus methods to more naturally occurring applications.  It is of no doubt that further applications of the fractional 
derivative and integral are yet to be discovered. In the future we plan to use such applications to further enrich the theory and mathematical power of fractional calculus.


\section*{Appendix} \label{Appendix}

It is left to show that (\ref{it_seq}) converges to (\ref{b series}), and that $b$ is actually twice differentiable with respect to time.  To prove this we show that each one of the series
\begin{align}
&\sum_{n=0}^\infty K_2^n(K_1)f \label{series}\\
&\sum_{n=0}^\infty\frac{\p}{\p t} K_2^n(K_1)f\label{d series}\\
&\sum_{n=0}^\infty\frac{\p^2}{\p t^2} K_2^n(K_1)f\label{d2 series}
\end{align}
converge uniformly, and to do so we will apply the Weierstrass $M$-test.

To apply Weierstrass we will first need to show that the function $
K_1f(x,t)
$
is continuous in $x$ and $t$.  We will actually show that it is twice differentiable with respect to $t$.  Now by definition
\begin{align*}
K_1(f(x,t))&=Ig(-\phi,-y)Ig(\phi,y)y^{-1}f\\
&=\int_0^t\underbrace{ g(-\phi,-y)\int_0^\tau \underbrace{g(\phi,y)y^{-1}f(x,s)}_1\ d s \ d\tau.}_2
\end{align*}
Since $f$ is continuous in $x$, so is $K_1f$.  Also since $f$ is  a continuous function of $t$, $g(\phi,y)$ is a composition of differentiable functions, and $y^{-1}$ is differentiable (this follows from the differentiability of $y$), then
\[
\int_0^\tau g(\phi,y)y^{-1}f(x,s)\ d s
\]
is a differentiable function of $\tau.$ For the same reasons the term labeled 2 is also the product of two differentiable functions with respect to time, thus 
\[
\int_0^tg(-\phi,-y)\int_0^\tau g(\phi,y)y^{-1}f(x,s)\ d s \ d\tau
\]
is differentiable with respect to $t$.  Therefore $K_1f$ is twice differentiable with respect to time, and continuous with respect to space. Now since $K_1f$ is a continuous function on a compact set, it is uniformly continuous and bounded.  Thus there exists a constant $C_1$ such that  for $(x,t)\in \mathcal{R}:= [0,1]\times [0,T]$, 
$
|K_1f|\le C_1
$.

We are now in a position to apply the $M$-test; first we show uniform convergence of (\ref{series}), and then move on to showing uniform convergence of (\ref{d series}) and (\ref{d2 series}).  As a first step towards showing (\ref{series}) converges uniformly, let 
\begin{align*}
C_2=\max_{\mathcal{R}}{\left|g(-\phi,-y),g(\phi, y)\right|} \quad \mbox{and} \quad C_3=\max_{\mathcal{R}}\left|{\frac{\p (gy^{-1})}{\p t}}y^{-1}\right|,
\end{align*}
 and observe:
\begin{align*}
K_2|(K_1f)| &\le  K_2 C_1\\
&\le  (Ig(-\phi,-y)J^\alpha g(\phi,y)+Ig(\phi,-y)J^\alpha I \frac{\p (ty^{-1})}{\p t}y)C_1\\
& = \frac{x^\alpha}{\Gamma(1+\alpha)}\left(tC_2^2+\frac{C_2C_3t^2}{2}\right)C_1 \\
&\le \frac{x^\alpha}{\Gamma(1+\alpha)}\left(TC_2^2+\frac{C_2C_3T^2}{2}\right)C_1.
\end{align*}
Therefore inductively
\begin{align*}
K_2^n|K_1|&\le  \frac{x^{n\alpha}}{\Gamma(1+n\alpha)}\left(TC_2^2+\frac{C_2C_3T^2}{2}\right)^nC_1\\
&\le  \frac{1}{\Gamma(1+n\alpha)}\left(TC_2^2+\frac{C_2C_3T^2}{2}\right)^nC_1
\end{align*}
and since
\[
\sum_{n=0}^\infty \frac{1}{\Gamma(1+n\alpha)}\left(TC_2^2+\frac{C_2C_3T^2}{2}\right)^nC_1
\]
is finite, our series converges uniformly.  We now show (\ref{d series}) is convergent.  The proof of this is quite similar.   Letting
\[
M_n= \frac{1}{\Gamma(1+n\alpha)}\left(TC_2^2+\frac{C_2C_3T^2}{2}\right)^n,
\]
one may show
\begin{align*}
\frac{\p}{\p t}K_2^n(K_1f)=&\frac{\p}{\p t}K_2(K_2^{n-1}(K_1f))\\
\le &\frac{\p}{\p t}K_2|(K_2^{n-1}(K_1f))|\\
\le &\frac{\p}{\p t} K_2 M_{n-1}C_1\\
\le &\Lambda M_{n-1}C_1.
\end{align*}
Thus (\ref{d series}) is also bounded above by a convergent series, so (\ref{d series}) is convergent.  An analogous argument may be given to show that the series (\ref{d2 series})  is convergent.

\begin{acknowledgement}
The first and third authors were partially supported by the National Science Foundation (USA) research grant DMS-1312529. The second author was partially supported by the U.S. Army Research Office grant W911NF-15-1-0537. 
\end{acknowledgement}

%
%
%
\bibliographystyle{elsarticle-num}
%
%
%
\bibliography{surfacevol}

\begin{thebibliography}{10}
\expandafter\ifx\csname url\endcsname\relax
  \def\url#1{\texttt{#1}}\fi
\expandafter\ifx\csname urlprefix\endcsname\relax\def\urlprefix{URL }\fi
\expandafter\ifx\csname href\endcsname\relax
  \def\href#1#2{#2} \def\path#1{#1}\fi

\bibitem{Herr}
R.~Herrmann, Fractional Calculus: An Introduction for Physicists, 2nd Edition,
  World Scientific, River Edge, NJ, 2014.

\bibitem{what}
J.~A. Machado, And {I} say to myself: ``{W}hat a fractional world!", Frac.
  Calc. Appl. Anal. 14~(4) (2011) 635--654.

\bibitem{quo}
J.~T. Machado, F.~Mainardi, V.~Kiryakova, Fractional calculus: Quo vadimus?
  ({W}here are we going?), Frac. Calc. Appl. Anal. 18~(2) (2015) 495--526.

\bibitem{kumar}
D.~Kumar, J.~Singh, S.~Kumar, A fractional model of {N}avier--{S}tokes equation
  arising in unsteady flow of a viscous fluid, Journal of the Association of
  Arab Universities for Basic and Applied Sciences 17 (2015) 14--19.

\bibitem{yang}
P.~Yang, Y.~C. Lama, Q.~Zhub, Constitutive equation with fractional derivatives
  for the generalized {UCM} model, J. Non-Newtonian Fluid Mech. 165 (2010)
  88--97.

\bibitem{Abel}
N.~H. Abel, Solution de quelques probl{\` emes \` a l'aide d'int\' egrales d\'
  efinies}, Mag. Naturv. 1~(2) (1823) 1--27.

\bibitem{Ross}
K.~S. Miller, B.~Ross, An Introduction to the Fractional Calculus and
  Fractional Differential Equations, Wiley, New York, 1993.

\bibitem{mainardi1}
F.~Mainardi, Y.~Luchko, G.~Pagnini, The fundamental solution of the space-time
  fractional diffusion equation, Frac. Calc. Appl. Anal. 4~(2) (2001) 153--192.

\bibitem{sie}
D.~Sierociuk, T.~Skovranek, M.~Macias, I.~Podlubny, I.~Petras, A.~Dzielinski,
  P.~Ziubinski, Diffusion process modeling by using fractional-order models,
  Appl. Math. Comput. 257~(15) (2015) 2--11.

\bibitem{tor}
P.~J. Torvik, R.~L. Bagley, On the appearance of the fractional derivative in
  the behavior of real materials, J. Appl. Mech. 51~(2) (1984) 294--298.

\bibitem{Caputo}
M.~Caputo, Linear model of dissipation whose {$Q$} is almost frequency
  independent - {II}, Geophys. J. R. Astr. Soc. 13 (1967) 529--539.

\bibitem{udita3}
U.~N. Katugampola, Mellin transforms of generalized fractional integrals and
  derivatives, Appl. Math. Comput. 257 (2015) 566--580.

\bibitem{oeis}
N.~J.~A. Sloane.
\newblock \href{http://oeis.org}{The on-line encyclopedia of integer sequences}
  [online] (2015).

\bibitem{Chen1}
H.~Chen, U.~N. Katugampola, Hermite–hadamard and hermite–hadamard–fejér
  type inequalities for generalized fractional integrals, J. Math. Anal. Appl.
  446~(2) (2017) 1274--1291.

\bibitem{u-1}
S.~Gaboury, R.~Tremblay, B.~{Fug\` ere}, Some relations involving a generalized
  fractional derivative operator, J. Inequal. Appl. (2013) 167.

\bibitem{mal}
A.~B. Malinowska, T.~Odzijewicz, D.~F.~M. Torres, Advanced Methods in the
  Fractional Calculus of Variations, Springer, New York, 2015.

\bibitem{u-3}
T.~Odzijewicz, A.~Malinowska, D.~Torres, A generalized fractional calculus of
  variations, Control and Cybernetics 42~(2) (2013) 443--458.

\bibitem{u-4}
T.~Odzijewicz, A.~Malinowska, D.~Torres, Fractional calculus of variations in
  terms of a generalized fractional integral with applications to physics,
  Abstract and Applied Analysis.

\bibitem{u-5}
A.~G. Butkovskii, S.~S. Postnov, E.~A. Postnova, Fractional integrodifferential
  calculus and its control-theoretical applications {I - M}athematical
  fundamentals and the problem of interpretation, Automation and Remote Control
  74~(4) (2013) 543--574.

\bibitem{u-5-2}
A.~G. Butkovskii, S.~S. Postnov, E.~A. Postnova, Fractional integrodifferential
  calculus and its control-theoretical applications. {II. F}ractional dynamic
  systems: Modeling and hardware implementation, Automation and Remote Control
  74~(5) (2013) 725--749.

\bibitem{Almeida1-doi}
R.~Almeida, \href{https://arxiv.org/abs/1601.07376}{Variational problems
  involving a caputo-type fractional derivative}, Journal of Optimization
  Theory and Applications\href {http://dx.doi.org/10.1007/s10957-016-0883-4}
  {\path{doi:10.1007/s10957-016-0883-4}}.
\newline\urlprefix\url{https://arxiv.org/abs/1601.07376}

\bibitem{mark}
R.~J. Marks, M.~W. Hall, Differintegral interpolation from a bandlimited
  signal's samples, IEEE Trans. Acoust. Speech Signal Process. 29 (1981)
  872--877.

\bibitem{bai}
J.~Bai, X.~C. Feng, Fractional-order anisotropic diffusion for image denoising,
  IEEE Trans. Image Process 16 (2007) 2492--2502.

\bibitem{ben}
D.~A. Benson, The fractional advection-dispersion equation{:} development and
  application (1998).

\bibitem{fre}
A.~D. Freed, K.~Diethelm, Y.~Luchko, Fractional-order viscoelasticity ({FOV}):
  {C}onstitutive development using the fractional calculus (first annual).,
  Technical Memorandum 2002-211914, Cleveland (2002).

\bibitem{mag}
R.~Magin, Fractional Calculus in Bioengineering, Begell House, Redding, 2006.

\bibitem{gloc}
W.~G. Gl\"{o}ckle, T.~F. Nonnenmacher, A fractional calculus approach to
  self-similar protein dynamics, Biophys. J. 68 (1995) 46--53.

\bibitem{gro1}
R.~Gorenflo, G.~D. Fabritiis, F.~Mainardi, Discrete random walk models for
  symmetric {L\' evy-F}eller diffusion processes, Physica A 269 (1999) 79--89.

\bibitem{gro2}
E.~Scalas, R.~Gorenflo, F.~Mainardi, Fractional calculus and continuous-time
  finance, Physica A 284 (2000) 376--384.

\bibitem{led}
C.~Lederman, J.-M. Roquejoffre, N.~Wolanski, Mathematical justification of a
  nonlinear integrodifferential equation for the propagation of spherical
  flames, C. R. Math. Acad. Sci. Paris 334 (2002) 569--574.

\bibitem{podl}
I.~Podlubny, Fractional-order systems and fractional-order controllers, Tech.
  Rep. UEF-03-94, Institute for Experimental Physics, Slovak Acad. Sci. (1994).

\bibitem{podl2}
I.~Podlubny, L.~Dorcak, J.~Misanek, Application of fractional-order derivatives
  to calculation of heat load intensity change in blast furnace walls, Trans.
  Tech. Univ. Ko$\check{s}$ice 5 (1995) 137--144.

\bibitem{jar}
A.~A. Jarbouh, Rheological behaviour modelling of viscoelastic materials by
  using fractional model, Energy Procedia 19 (2012) 143--157.

\bibitem{ahmad}
W.~M. Ahmad, R.~El-Khazali, Fractional-order dynamical models of love, Chaos
  Solitons Fractals 33 (2007) 1367--1375.

\bibitem{song}
L.~Song, S.~Y. Xu, J.~Y. Yang, Dynamical models of happiness with fractional
  order, Commun. Nonlinear Sci. Numer. Simulat. 15 (2010) 616--628.

\bibitem{Edwards2}
D.~A. Edwards, Transport effects on surface reaction arrays: Biosensor
  applications, Mathematical Biosciences 230 (2011) 12--22.

\bibitem{pod}
I.~Podlubny, Fractional Differential Equations: Mathematics In Science and
  Engineering, Academic Press, San Diego, 1999.

\bibitem{samko}
S.~G. Samko, A.~A. Kilbas, O.~I. Marichev, Fractional Integrals and
  Derivatives. Theory and Applications, Gordon and Breach, Amsterdam, 1993.

\bibitem{udita1}
U.~N. Katugampola, New approach to a generalized fractional integral, Appl.
  Math. Comput. 218~(3) (2011) 860--865.

\bibitem{udita2}
U.~N. Katugampola, New approach to generalized fractional derivatives, Bull.
  Math. Anal. Appl. 6~(4) (2014) 1--15.

\bibitem{kol1}
K.~M. Kolwankar, Local fractional calculus: A review, arXiv:1307.0739. (2013).

\bibitem{kol2}
K.~M. Kolwankar, A.~D. Gangal, Fractional differentiability of nowhere
  differentiable functions and dimensions, Chaos 6~(4) (1996) 505--513.

\bibitem{katu}
D.~R. Anderson, D.~J. Ulness, Properties of the {K}atugampola fractional
  derivative with potential application in quantum mechanics, J. Math. Phys 56
  (2015) 063502.

\bibitem{wfra}
M.~D. Ortigueiraa, J.~A. Machado, What is a fractional derivative?, J.
  Computational Physics 293 (2015) 4--13.

\bibitem{udita4}
U.~N. Katugampola, Correction to “what is a fractional derivative?” by
  ortigueira and machado [j. comput. phys. 293(2015):4–13. special issue on
  fractional pdes], J. Comput. Phys. 321 (2016) 1255--1257.

\bibitem{kai}
K.~Diethelm, The Analysis of Fractional Differential Equations{:} An
  Application-Oriented Exposition Using Differential Operators of Caputo Type,
  Springer, New York, 2010.

\bibitem{pod1}
I.~Podlubny, Geometric and physical interpretation of fractional integration
  and fractional differentiation, Frac. Calc. Appl. Anal. 5~(4) (2002)
  367--386.

\bibitem{pod2}
I.~Podlubny, V.~Despotovic, T.~Skovranek, B.~H. McNaughton, Shadows on the
  walls: Geometric interpretation of fractional integration, Journal of Online
  Mathematics and Its Applications 7, article ID 1664.

\bibitem{adda1}
F.~{Ben Adda}, Geometric interpretation of the fractional derivative, Journal
  of Fractional Calculus 11 (1997) 21--52.

\bibitem{adda2}
F.~{Ben Adda}, Interpretation geometrique de la differentiabilite et du
  gradient d'ordre reel, C. R. Acad. Sci. Paris 326~(I) (1998) 931--934.

\bibitem{Ren}
F.~Ren, Z.~Yu, F.~Su., Fractional integral associated to the self-similar set
  of the generalized self-similar set and its physical interpretation, Phys.
  Lett. A 219 (1996) 59--68.

\bibitem{Gor1}
R.~Gorenflo, Afterthoughts on interpretation of fractional derivatives and
  integrals., in: P.~Rusev, I.~Dimovski, V.~Kiryakova (Eds.), Transform Methods
  and Special Functions Varna'96 (Proc. 3rd Internat. Workshop), Bulgarian
  Academy of Sciences, Institute of Mathematics and Informatics, Sofia, 1998,
  pp. 589--591.

\bibitem{Hey}
N.~Heymans, I.~Podlubny, Physical interpretation of initial conditions for
  fractional differential equations with {R}iemann--{L}iouville fractional
  derivatives, Rheologica Acta 45~(5) (2006) 765--772.

\bibitem{Kirk1}
V.~Kiryakova, A long standing conjecture failed?, in: P.~Rusev, I.~Dimovski,
  V.~Kiryakova (Eds.), Transform Methods and Special Functions Varna'96 (Proc.
  3rd Internat. Workshop), Institute of Mathematics and Informatics, Bulgarian
  Academy of Sciences, Sofia, 1998, pp. 579--588.

\bibitem{Mon}
M.~Monsrefi-Torbati, J.~K. Hammond, Physical and geometrical interpretation of
  fractional operators, J. Franklin Inst. 335B~(6) (1998) 1077--1086.

\bibitem{Nig}
R.~R. Nigmatullin, A fractional integral and its physical interpretation,
  Theoret. and Math. Phys. 90~(3) (1992) 242--251.

\bibitem{Rut}
R.~S. Rutman, On the paper by {R. R. N}igmatullin `{A} fractional integral and
  its physical interpretation', Theoret. and Math. Phys. 100~(3) (1994)
  1154--1156.

\bibitem{Rut2}
R.~S. Rutman, On physical interpretations of fractional integration and
  differentiation, Theoret. and Math. Phys. 105~(3) (1995) 1509--1519.

\bibitem{Yu1}
Z.~Yu, F.~Ren, J.~Zhou, Fractional integral associated to generalized
  cookie-cutter set and its physical interpretation, J. Phys. A: Math. Gen. 30
  (1997) 5569--5577, 

\bibitem{Mac}
J.~A. Machado, A probabilistic interpretation of the fractional order
  differentiation, Frac. Calc. Appl. Anal. 6~(1) (2003) 73--80.

\bibitem{recent1}
M.~H. Tavassoli, A.~Tavassoli, M.~R. {Ostad Rahimi}, The geometric and physical
  interpretation of fractional order derivatives of polynomial functions,
  Differential Geometry - Dynamical Systems 15 (2013) 93--104.

\bibitem{recent2}
S.~T. Nizami, N.~Khan, F.~H. Khan, A new approach to represent the geometric
  and physical interpretation of fractional order derivatives of polynomial
  function and its application in field of sciences, Canadian Journal on
  Computing in Mathematics, Natural Science, Engineering \& Medicine 1~(1)
  (2010) 1--8.

\bibitem{herr2}
R.~Herrmann, Towards a geometric interpretation of generalized fractional
  integrals---{Erd\' e}lyi-{K}ober type integrals on {$R^N$}, as an example,
  Frac. Calc. Appl. Anal. 17~(2) (2014) 361--370.

\bibitem{Main2}
F.~Mainardi, Considerations on fractional calculus: Interpretations and
  applications, in: P.~Rusev, I.~Dimovski, V.~Kiryakova (Eds.), Transform
  Methods and Special Functions Varna'96 (Proc. 3rd Internat. Workshop),
  Institute of Mathematics and Informatics, Bulgarian Academy of Sciences,
  Sofia, 1998, pp. 594--597.

\bibitem{jake}
A.~J. Jakeman, R.~S. Anderssen, Abel type integral equations in stereology {I.
  G}eneral discussion 105~(2) (1975) 121--133.

\bibitem{kow1}
G.~Kowalewski, Integralgleichungen, Water de Gruyter and Co, Berlin, 1930.

\bibitem{avaz}
Z.~Avazzadeh, B.~Shafiee, G.~B. Loghmani, Fractional calculus of solving
  {A}bel's integral equations using {C}hebyshev polynomials, Appl. Math. Sci
  5~(45) (2011) 2207--2216.

\bibitem{Mine}
G.~N. Minerbo, M.~E. Levy, Inversion of {A}bel's integral equation by means of
  orthogonal polynomials, SIAM J. Numer. Anal. 6~(4) (1969) 598--616.

\bibitem{paul}
P.~P.~B. Eggermont, On galerkin methods for abel-type integral equations, SIAM
  J. Numer. Anal. 25~(5) (1988) 1093--1117.

\bibitem{brun}
H.~Brunner, Collocation Methods for Volterra Integral and Related Functional
  Differential Equations, Vol.~15 of Cambridge Monographs on Applied and
  Computational Mathematics, Cambridge Univ Press, Cambridge, 2004.

\bibitem{lep}
U.~Lepik, Solving fractional integral equations by the {H}aar wavelet method,
  Appl. Math. Comput. 214~(2) (2009) 468--478.

\bibitem{sae}
H.~Saeedi, N.~Mollahasani, M.~M. Moghadam, G.~N. Chuev, An operational {H}aar
  wavelet method for solving fractional volterra integral equations, Int. J.
  Appl. Math. Comput. Sci. 21~(3) (2011) 535--547.

\bibitem{sae2}
H.~Saeedi, M.~M. Moghadam, N.~Mollahasani, G.~N. Chuev, A {CAS} wavelet method
  for solving nonlinear {F}redholm integrodifferential equations of fractional
  order, Commun. Nonlinear Sci. Nummer. Simul. 16~(3) (2011) 1154--1163.

\bibitem{Li}
M.~Li, W.~Zhao, Solving {A}bel's type integral equation with {M}ikusinski's
  operator of fractional order, Advances in Mathematical PhysicsArticle ID
  806984.

\bibitem{Kanwal}
R.~P. Kanwal, Linear Integral Equations{:} Theory \& Technique, Springer, New
  York, 2013.

\bibitem{Salman}
S.~Jahanshahi, E.~Babolian, D.~F.~M. Torres, A.~Vahidi, Solving {A}bel integral
  equations of first kind via fractional calculus, J. King Saud University -
  Science 27 (2015) 161--167.

\bibitem{Gal}
I.~M. Gelfand, K.~Vilenkin, Generalized Functions, Vol.~1, Academic Press, New
  York, 1964.

\bibitem{Edwards1}
D.~A. Edwards, Transport effects on surface-volume biological reactions,
  Journal of Mathematical Biology 39~(6) (1999) 533--561.

\bibitem{Edwards3}
D.~A. Edwards, Testing the validity of the effective rate constant
  approximation for surface reaction with transport, Applied Mathematics
  Letters 15 (2002) 547--552.

\bibitem{Loverro}
A.~Loverro, Fractional Calculus: History, Definitions, and Applications for the
  Engineer, 2004, rapport Technique.

\bibitem{brun2}
H.~Brunner, A survey of recent advances in the numerical treatment of volterra
  integral and integrodifferential equations, J. Comput. Appl. Math. 8 (1982)
  213--229.

\bibitem{olm}
W.~E. Olmstead, R.~A. Handelsman, Diffusion in a semi-infinite region with
  nonlinear surface dissipation, SIAM Rev. 18 (1976) 275--291.

\bibitem{Grabowski}
E.~Grabowski, L.~Friedman, E.~Leonard, Effects of shear rate on diffusion and
  adhesion of blood platelets to a foreign surface, Ind. Eng. Chem. Fund. 11
  (1972) 224--232.

\bibitem{Bertucci}
C.~Bertucci, A.~Piccoli, M.~Pistolozzi, Optical biosensors as a tool for early
  determination of absorption of lead candidates and drugs, Combinatorial
  Chemistry and High Throughput Screening 10 (2007) 433--440.

\bibitem{Raghaven}
M.~Raghavan, M.~Y. Chen, L.~N. Gastinel, P.~J. Bjorkman, Investigation of the
  interaction between the class {I} {MHC}-related {F}c receptor and its
  immunoglobulin {G} ligand, Immunity 1 (1994) 303--315.

\bibitem{Rich}
R.~L. Rich, D.~G. Myszka, Survey of the year 2009 commercial optical biosensor
  literature, J. Mol. Recognit. 24 (2011) 892--914.

\end{thebibliography}
\end{document}